**Transmuted Geometric Distribution and its Prpoerties**

Subrata Chakraborty

Department of Statistics, Dibrugarh University, Dibrugarh-786004, Assam, India

E-mail address: subrata_arya@yahoo.co.in

**Abstract** *Transmuted geometric distribution with two parameters $q\,(0 < q < 1)$ and $\alpha\,(> 0)$ is proposed as a new generalization of the geometric distribution by employing the quadratic transmutation techniques of Shaw and Buckley (2007). Its important distributional and reliability properties are investigated. Parameter estimation methods is discussed.*

**1**. **Introduction**

A random variable $X$ follows geometric distribution with parameter $q$ [GD($q$)] (Johnson et al., 2005, page 210, equation (5.8)) if its probability mass function(pmf) is given by

$$P(X = t) = p\,q^t,\ t = 0,1,2,\cdots; 0 < q < 1; p = 1 - q \qquad (1)$$

For geomtric distribution in (1) the following hold.

i. cumulative distribution function(cdf): $F_X(t) = \Pr(X \le t) = 1 - q^{t+1}$

ii. survival function: $S_X(t) = \Pr(X \ge t) = q^t$

iii. hazard rate function: $r_X(t) = \Pr(X = t)/S_X(t) = p$

iv. reversed hazard rate function: $r_X^*(t) = P(X = t)/F_X(t) = p\,q^t/(1 - q^{t+1})$

Many generalizations of geometric distribution has been attempted by researchers by using different methods (see Chakraborty and Gupta 2013 among others).

Transmutation method first introduced by Shaw and Buckley (2007) has been used by many authors to generate new distrbutions staring with suitable continuous distributions (see Oguntunde and Adejumo, 2015 for details). But sofar there is no evidence of any attempt to use this method to generate new discrete distribution. In this article an attempt has been made to derive a new generalization of geometric distrbution with two parameters $0 < q < 1$ and $-1 < \alpha < 1$ by using the quadratic rank transmutation map method of Shaw and Buckley (2007). Some important distributional and reliability properties, classification of hazard rate, a derivation as discrete concentartion, parameter estimation methods have been discussed.

## 2. Transmuted geometric distribution

A random variable $Y$ is said to be constructed by the quadratic rank transmutation map method of Shaw and Buckley (2007) by transmuting another random variable $X$ with cdf $F_X()$ if the cumulative distribution function (cdf) of $Y$ is given by

$$F_Y(y) = (1+\alpha)F_X(y) - \alpha(F_X(y))^2, -1 < \alpha < 1$$

Suppose the random variable $X$ has geometric distribution (Johnson et al., 2005, page 210, equation (5.8)) with pmf and cdf given in section 1.

Then the cdf of the transmuted geometric variable $Y$ will be constructed as

$$F_Y(y) = (1+\alpha)(1-q^{y+1}) - \alpha(1-q^{y+1})^2, y = 1, 2, \cdots; 0 < q < 1; -1 < \alpha < 1$$
$$= 1 + (\alpha-1)q^{y+1} - \alpha q^{2(y+1)}.$$

The corresponding pmf will then be given by

$$P(Y = y) = (1-\alpha)q^y(1-q) + \alpha q^{2y}(1-q^2), y = 1, 2, \cdots; 0 < q < 1; -1 < \alpha < 1 \quad (2)$$

The distribution in (2) will hence forth be refered to as the Transmuted geometric distribution(TGD) with two parameters $q$ and $\alpha$. In short $TGD(q, \alpha)$.

**Particular cases**:

(i) $\alpha = -1$, (2) reduces to a special case of the Exponentiated Geometric distribution $(q, 2)$ (Chakraborty and Gupta, 2013) which is the *distribution of the maximum* of two independent $GD(q)$c random variable.

(ii) $\alpha = 0$, (2) reduces to $GD(q)$ in (1)

(iii) $\alpha = 1$, to $GD(q^2)$ with pmf $q^{2y}(1-q^2)$ which is the *distribution of the minimum* of two independent $GD(q)$ random variable.

**Remark 1.** $TGD(q, \alpha)$ forms a continuous bridge between the distribution of the minimum to maximum of sample of size two from geometric distribution.

### 2.2. Distributional Properties

#### 2.2.1. Shape of the pmf

The graphs of the pmf are plotted for various combinations of the values of the two parameters $q$ and $\alpha$ in Figure 1.

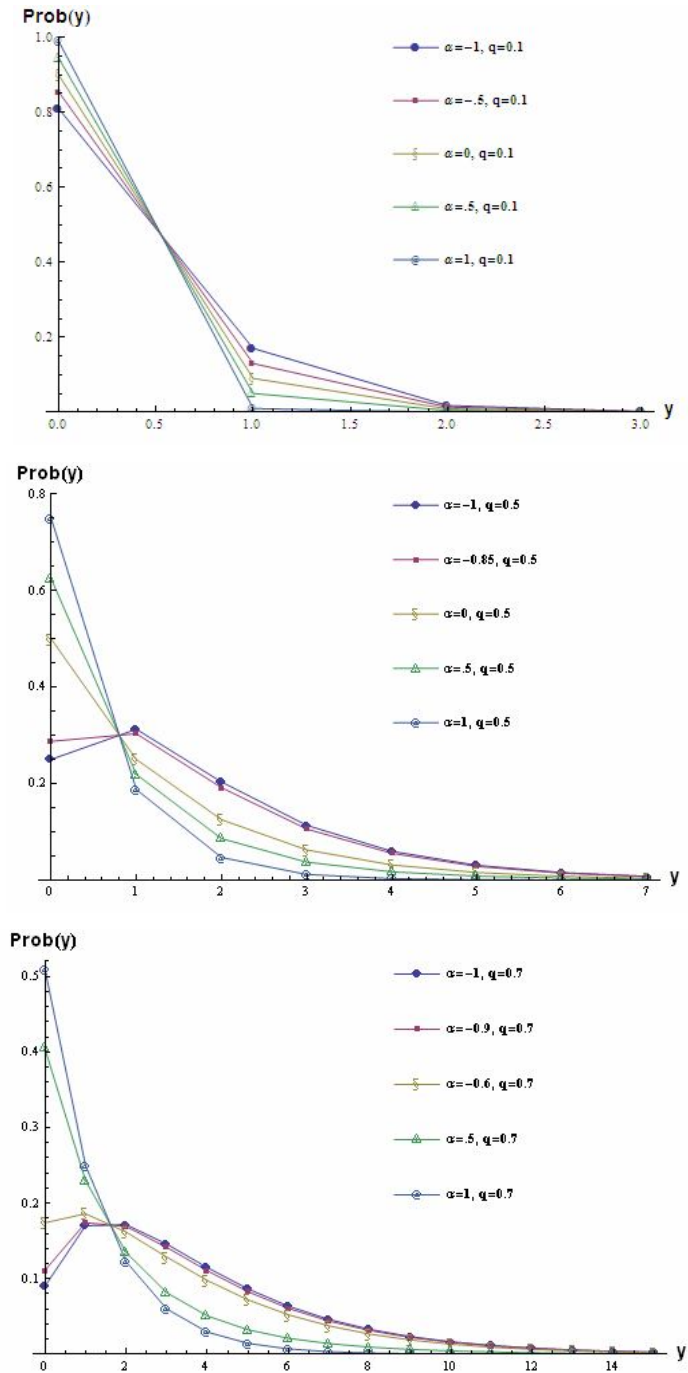

**Figure 1. pmf plot of TGD**$(q,\alpha)$

When $-1 \leq \alpha \leq 0$, the pmf is unimodal if $-1 < \alpha < -\{q(2+q)\}^{-1}$ and decreasing if $-\{q(2+q)\}^{-1} < \alpha < 1$ provided $q > 0.414$. For $0 \leq \alpha \leq 1$ the pmf is always a decreasing function. The mode shifts form lower to higher values with increasing values of the parameter $q$ and $\alpha$. With increase in the value of $q$ and $\alpha$ the spread of the distribution also increases. The above assertions are mathematically established in theorem 1.

**2.2.2. Mode**

**Theorem 1.** $\text{TGD}(q,\alpha)$ is unimodal for $-1 < \alpha < -\{q(2+q)\}^{-1}$ provided $q > 0.414$.

*Proof*: A pmf $P(Y = y)$ with support $y = 0, 1, 2, \cdots$, is unimodal if there exist a unique point $M (\neq 0)$, in the support of $Y$ such that $P(Y = y)$ is increasing on $(0, 1, \cdots, M)$ and decreasing on $(M, M+1, \cdots)$. $M$ is then the unique mode of $P(Y = y)$.

Thus to be unimodal we must have

$P(Y = 1) > P(Y = 0)$

$\Rightarrow (1-\alpha)q(1-q) + \alpha q^2(1-q^2) > (1-\alpha)(1-q) + \alpha(1-q^2)$

$\Rightarrow (1-\alpha)(1-q)^2 + \alpha(1-q^2)(1-q^2) < 0$

$\Rightarrow \alpha < -(1-q)^2 / \{(1-q^2)^2 - (1-q)^2\} = -1/\{(1+q)^2 - 1\} = -1/\{q(2+q)\}$

But the condition $-1 < \alpha < -\{q(2+q)\}^{-1}$ makes sense only if $q(2+q) > 1$ which implies $q > \sqrt{2} - 1 \cong 0.414$. Hence proved.

**Remark 2.** For $q < 0.414$ the condition of unimodality leads to $\alpha$ out side its permissible range of (-1, +1).

**Remark 3.** For $0 \leq \alpha \leq 1$, the pmf is decreasing with the mode occuring at the point '0'.

**2.2.3. *An alternative derivation of the pmf***

**Theorem 2.** $\text{TGD}(q,\alpha)$ is the discrete analogue of the skew exponential of Shaw and Buckley (2007).

*Proof*:

The pdf and cdf of the skew exponential distribution derived using the quadratic rank transmutation are respectively given by

$f_X(x) = \beta e^{-\beta x}(1-\alpha) + 2\alpha\beta e^{-2\beta x}, x > \beta > 0, -1 < \alpha < 1$

$F_X(x) = (1+\alpha)(1 - e^{-\beta x}) - \alpha(1 - e^{-\beta x})^2, x > \beta > 0, -1 < \alpha < 1$ (Shaw and Buckley, 2007)

Hence the pmf of the discrete analogue $Y = \lfloor X \rfloor$, where $\lfloor X \rfloor$ is the floor function, is given by the formula

$P(Y = y) = F_X(y+1) - F_X(y)$ [see Chakraborty and Charavarty, 2015 for details]. This on simplification reduces to the pmf of $\text{TGD}(q,\alpha)$ with $e^{-\beta} = q$.

**2.2.4 *Probability genrating function***

**Theorem 3.** The probability generating function of $\text{TGD}(q,\alpha)$ is given by

$$\frac{(1-q)(1-\alpha q(1-z)-q^2 z)}{(1-qz)(1-q^2 z)}, \ |qz|<1$$

*Proof*:

It is known that the pgf $E[z^X]$ of $GD(q)$ is equal to $\frac{(1-q)}{(1-qz)}$ (page 215, Johnson et al., 2005)

Therefore pgf of $TGD(q,\alpha)$ is given by

$$E(z^Y) = \sum_y z^y P(Y=y) = \sum_y z^y \{(1-\alpha)q^y(1-q) + \alpha q^{2y}(1-q^2)\}$$

$$= (1-\alpha)(1-q)\sum_y z^y q^y + \alpha(1-q^2)\sum_y (zq^2)^y$$

$$= \frac{(1-\alpha)(1-q)}{(1-qz)} + \frac{\alpha(1-q^2)}{(1-q^2 z)}, \ |qz|<1$$

The result follows on simplification.

### 2.2.5 *Moment and related measures*

**Theorem 4.** The $r^{th}$ factorial moment $TGD(q,\alpha)$ is given by

$$E[Y_{(r)}] = (1-\alpha)r!\left(\frac{q}{1-q}\right)^r + \alpha r!\left(\frac{q^2}{1-q^2}\right)^r.$$

*Proof*:

It is known that the $r^{th}$ factorial moment $E[X_{(r)}] = E[X(X-1)...(X-r+1)]$ of $GD(q)$ is given by (page 215, Johnson et al., 2005)

$$E(X_{(r)}) = r!\left(\frac{q}{1-q}\right)^r \tag{3}$$

Therefore $r^{th}$ factorial moment of $TGD(q,\alpha)$ is given by

$$E(Y_{(r)}) = \sum_y y^r P(Y=y) = \sum_y y^r \{(1-\alpha)q^y(1-q) + \alpha q^{2y}(1-q^2)\}$$

$$= (1-\alpha)\sum_y y^r q^y (1-q) + \alpha \sum_y y^r (q^2)^y (1-q^2)$$

$$= (1-\alpha)r!\left(\frac{q}{1-q}\right)^r + \alpha r!\left(\frac{q^2}{1-q^2}\right)^r \ [\text{using (3)}]$$

In particular

$$E(Y) = (1-\alpha)\left(\frac{q}{1-q}\right) + \alpha\left(\frac{q^2}{1-q^2}\right) = \frac{q(1-\alpha)+q^2}{1-q^2}$$

$$E(Y(Y-1)) = (1-\alpha)\frac{2q^2}{(1-q)^2} + \frac{2\alpha q^4}{(1-q^2)^2} = \frac{2q^2}{(1-q)^2}\left[(1-\alpha) + \frac{\alpha q^2}{(1+q)^2}\right]$$

$$= \frac{2q^2}{(1-q)^2}\left[\frac{(1-\alpha)(1+q)^2 + \alpha q^2}{(1+q)^2}\right] = \frac{2(1-\alpha)q^2 + 4(1-\alpha)q^3 + 2q^4}{(1-q^2)^2}$$

**Theorem 5.** The $r^{th}$ factorial cumulant $TGD(q,\alpha)$ is given by

$$(1-\alpha)(r-1)!\left(\frac{q}{1-q}\right)^r + \alpha(r-1)!\left(\frac{q^2}{1-q^2}\right)^r.$$

**Various moments from pgf:**

**Factorial moments :**

$$E[Y_{(2)}] = \frac{2q^2\{(1+q)^2 - \alpha(1+2q)\}}{(1-q^2)^2}$$

$$E[Y_{(3)}] = 6q^3\left[\frac{(1-\alpha)}{(1-q^2)} + \frac{\alpha q^3}{(1-q^2)^3}\right]$$

$$E[Y_{(4)}] = 24q^4\left[\frac{\alpha q^4}{(1-q^2)^4} - \frac{1-\alpha}{(1-q)^4}\right]$$

**Raw moments:**

$$E[Y^2] = \frac{q[(1+q)^3 - \alpha\{q(3q+2)+1\}]}{(1-q^2)^2}$$

$$E[Y^3] = \frac{q(1+q)^3\{q(q+4)+1\} - \alpha q[q\{q(q(7q+12)+16)+6\}+1]}{(1-q^2)^3}$$

$$E[Y^4] = \frac{q(1+q)^5\{q(q+10)+1\} - \alpha q[q(q(q(5q(q(3q+10)+23)+104)+61)+14)+1)]}{(1-q^2)^4}$$

**Central moments (moments about mean):**

$$E[(Y-E(Y))^2] = \frac{q\{1-\alpha^2 + q(1-\alpha^2 + q(1-\alpha)+2)\}}{(1-q^2)^2}$$

$$E[(Y-E(Y))^3] = -\frac{q\{2\alpha^3 q^2 + 3\alpha^2 q(1+q^2) + \alpha(q^4 + 4q^2 + 1) - (1+q)^4\}}{(1-q^2)^3}$$

$$E[(Y-E(Y))^4] = \{q(\lambda^4 q^3 - 2\lambda^3 q^2(q(11q+8)+3) + 2\lambda^2 q(q(q(q(38q+69)+62)+21)+2)$$
$$- \lambda(q(q(q(q(q(81q+262)+423)+332)+135)+22)+1)$$
$$+ (3q(q(4q+11)+5)+1)(q+1)^2\}/(q^2-1)^4$$

**Index of dispersion** :

$$\text{ID} = \text{Varince / Mean} = \frac{q\{1 - \alpha^2 + q(1 - \alpha^2 + q(1 - \alpha) + 2)\}}{(1 - q^2)\{q(1 - \alpha) + q^2\}}$$

**Remark 4.** It can be easily seen that the proposed distribution is always over dispersed (ID greater than 1).

**Skeweness and Kurtosis:**

$$\beta_1 = \frac{[q\{2\alpha^3 q^2 + 3\alpha^2 q(1 + q^2) + \alpha(q^4 + 4q^2 + 1) - (1 + q)^4\}]^2}{[q\{1 - \alpha^2 + q(1 - \alpha^2 + q(1 - \alpha) + 2)\}]^3}$$

$$\beta_2 = [\{q(\lambda^4 q^3 - 2\lambda^3 q^2 (q(11q + 8) + 3) + 2\lambda^2 q(q(q(q(38q + 69) + 62) + 21 + 2)$$
$$- \lambda(q(q(q(q(q(81q + 262) + 423) + 332) + 135) + 22) + 1)$$
$$+ (3q(q(4q + 11) + 5) + 1)(q + 1)^2\}] / q\{1 - \alpha^2 + q(1 - \alpha^2 + q(1 - \alpha) + 2)\}$$

*2.2.6 Median*

**Theorem 6.** The median of $\text{TGD}(q, \alpha)$ is given by

$$\left\lfloor \frac{\log\left\{(\alpha - 1) + \sqrt{(\alpha^2 + 1)}\right\} - \log\{2\alpha\}}{\log q} \right\rfloor, \text{ where } \lfloor y \rfloor \text{ is the floor function which the largest}$$

integer less than $y$.

*Proof* :

The median can be obtained by solving the following equation

$$\alpha q^{2(y_{.5} + 1)} - (\alpha - 1)q^{y_{.5} + 1} = 0.5$$

$$\Rightarrow \{\alpha q^{y_{.5} + 1} - (\alpha - 1)\} q^{y_{.5} + 1} = 0.5$$

$$\Rightarrow \{\alpha z - (\alpha - 1)\} z = 0.5$$

$$\Rightarrow \alpha z^2 - (\alpha - 1)z - 0.5 = 0$$

$$\Rightarrow z = \frac{(\alpha - 1) \pm \sqrt{\alpha^2 + 1}}{2\alpha}$$

$$\therefore y_{.5} = \left\lfloor \frac{\log\left\{(\alpha - 1) + \sqrt{(\alpha^2 + 1)}\right\} - \log\{2\alpha\}}{\log q} \right\rfloor$$

[since $(\alpha - 1) - \sqrt{(\alpha^2 + 1)}$ leads to indetrminancy].

**Remark 5.** The $p^{th}$ quantile $y_p$ can similarly obtained as solution of the equation

$$\alpha q^{2(y_p + 1)} - (\alpha - 1)q^{y_p + 1} = 1 - p$$

## 2.3 Reliabilty properties

### 2.3.1. Survival function

$$S_Y(y) = P(Y \geq y) = (1-\alpha)q^y + \alpha q^{2y}$$

### 2.3.2. Hazard rate function and its classification

$$r_Y(y) = \Pr(Y = y)/S_Y(y)$$

$$= \frac{(1-\alpha)q^y(1-q) + \alpha q^{2y}(1-q^2)}{(1-\alpha)q^y + \alpha q^{2y}}$$

$$= \frac{(1-\alpha)(1-q) + \alpha q^y(1-q^2)}{(1-\alpha) + \alpha q^y}$$

**Theorem 7.** *The* $\text{TGD}(q,\alpha)$ *has increasing, decreasing and constant hazard rate for* $-1 < \alpha < 0$, $0 < \alpha < 1$ *and* $\alpha = 0 \text{ or } 1$ *respetively.*

*Proof* :

$$r_Y(y) = \frac{(1-\alpha)(1-q) + \alpha q^y(1-q^2)}{(1-\alpha) + \alpha q^y}$$

$$= \frac{(1-\alpha) + \alpha q^y - \{(1-\alpha)q + \alpha q^{y+2}\}}{(1-\alpha) + \alpha q^y}$$

$$= 1 - q\frac{(1-\alpha) + \alpha q^{y+1}}{(1-\alpha) + \alpha q^y}$$

But $q\dfrac{(1-\alpha) + \alpha q^{y+1}}{(1-\alpha) + \alpha q^y}$ is a decreasing (increasing) function of $y$ for $-1 < \alpha < 0$ ($0 < \alpha < 1$).

Hence $r_Y(y)$ is increasing (decreasing) function of $y$ for $-1 < \alpha < 0$ ($0 < \alpha < 1$).

Constant hazard rates are obtained as $r_Y(y) = q$ when $\alpha = 0$ and $r_Y(y) = q^2$ for $\alpha = 1$.

### 2.3.3 Plots of hazard rate function

The hazard rate function of $\text{TGD}(q,\alpha)$ are plotted in figure 4 for various values of the parameters to demonstrate the monotonic properties.

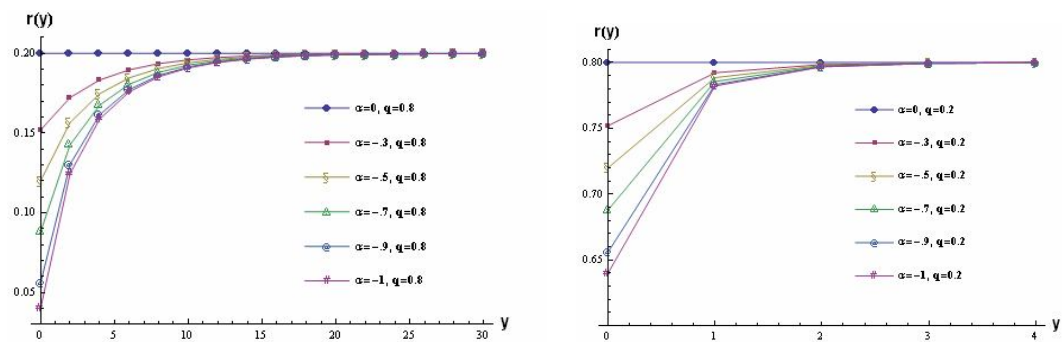

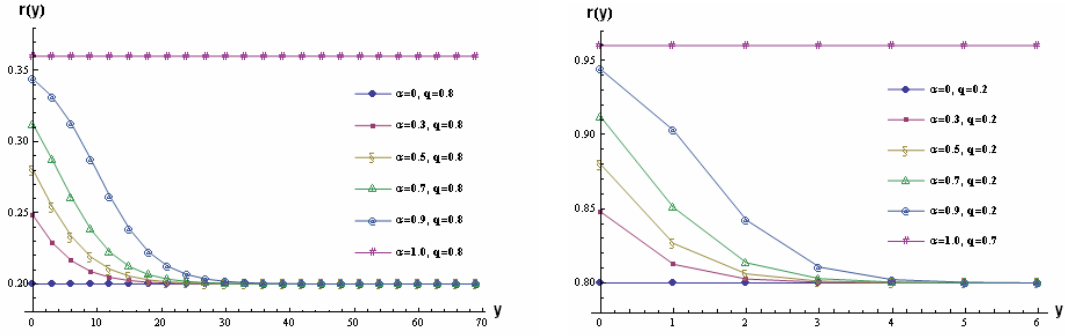

**Figure 2. Hazard rate Function plots TGD$(q, \alpha)$**

From the plots it is observed that the hazard rate of the TGD$(q,\alpha)$ is **increasing** for $-1 < \alpha < 0$, **decreasing** when $0 < \alpha < 1$ and **constant** if $\alpha = 0$ or $1$ Also it is seen that even when $\alpha \neq 1$, the hazard rate tend to approach constant as $t$ increses. Smaller the value of $q$ the faster is the rate of stabilization of the hazard rate.

### 2.3.3 Reversed hazard rate function

$$r_Y^*(y) = P(Y=y)/F_Y(y) = \frac{(1-\alpha)q^y(1-q) + \alpha q^{2y}(1-q^2)}{1-(1-\alpha)q^{y+1} - \alpha q^{2y+2}}$$

### 3. Parameter Estimation

**3.1 *From sample proportion of 1s and 0s*:** If $p_0, p_1$ be the known observed proportion of 0's and 1's in the sample then the parameters can be estimated solving the equations:

$$p_0 = (1-\alpha)(1-q) + \alpha(1-q^2) \text{ and } p_1 = (1-\alpha)q(1-q) + \alpha q^2(1-q^2)$$

**3.2. *From sample quantiles*:** If $t_1$, $t_2$ be two observed points such that $F_Y(t_1) = p_1, F_Y(t_2) = p_2$, then the two parameters $q$ and $\alpha$ can be estimated by solving the simultaneous equations

$$p_1 = 1 + (\alpha-1)q^{t_1+1} - \alpha q^{2(t_1+1)} \text{ and } p_2 = 1 + (\alpha-1)q^{t_2+1} - \alpha q^{2(t_2+1)}$$

### 3.3 *Method of moments*:

The moment estimates are obtained following the method proposed by Khan et al. (1989). Here the moment estimates of $q$ and $\alpha$ are obtained by minimizing $(E(Y) - m_1)^2 + (E(Y^2) - m_2)^2$ with respect to $q$ and $\alpha$, where $m_1$ and $m_1$ are the 1$^{st}$ and the 2$^{nd}$ observed raw moments respectively.

$$= \left(\frac{q(1-\alpha) + q^2}{1-q^2} - m_1\right)^2 + \left(\frac{q[(1+q)^3 - \alpha\{q(3q+2)+1\}]}{(1-q^2)^2} - m_2\right)^2$$

**3.4 Maximum likelihood estimation**:

For a sample $(y_1, y_2, \cdots, y_n)$ of size $n$ is drawn from $\text{TGD}(q,\alpha)$ the likelihood function is given by $L = \prod_{i=1}^{n} \{(1-\alpha)q^{y_i}(1-q) + \alpha q^{2y_i}(1-q^2)\}$. Taking logarithm on both sides gives

$$L = \prod_{i=1}^{n} q^{y_i}(1-q)\{(1-\alpha) + \alpha q^{y_i}(1+q)\}.$$

Therefore

$$\log L = \sum_{i=1}^{n} \log[q^{y_i}(1-q)\{(1-\alpha) + \alpha q^{y_i}(1+q)\}]$$

$$= \sum_{i=1}^{n} \log[q^{y_i}(1-q)] + \sum_{i=1}^{n} \log[(1-\alpha) + \alpha q^{y_i}(1+q)]$$

$$= n\log(1-q) + n\bar{y}\log[q] + \sum_{i=1}^{n} \log[(1-\alpha) + \alpha q^{y_i}(1+q)]$$

Since the likelihood equations have no closed form solution, ML estimates can be obtained by maximizing log likelihood function using global numerical maximization techniques.

**4. Conclusion**

A new generalization of the geometric distribution is proposed using quadratic rank transmuation and some of its prperties presented. Further works on data modeling, inference with real life applications will be communicated in follow up articles.